\documentclass[pre,preprintnumbers,amsmath,amssymb]{revtex4}
\usepackage{color,hyperref}
\usepackage{graphicx}
\usepackage{epsfig}
\usepackage{epstopdf}
\usepackage{bm}

\newtheorem{theorem}{Theorem}[section]

\newtheorem{lemma}{Lemma}[section]

\begin{document}

\title{Existence of exponentially spatially localised breather solutions for  lattices of nonlinearly coupled  particles: Schauder's fixed point theorem approach}

\author{Dirk Hennig}
\author{Nikos I. Karachalios}
\affiliation{Department of Mathematics, University of Thessaly, Lamia GR 35100, Greece}


\begin{abstract}
\noindent {\cal {\bf Abstract:}} The problem of showing the existence of localised modes in nonlinear lattices has attracted considerable efforts from the physical but also from the mathematical viewpoint where a rich variety of methods has been employed.  In this paper we prove that a fixed point theory approach based on the celebrated Schauder's Fixed Point Theorem may provide a general method to establish concisely not only the existence of localised structures but also  a required rate of spatial localisation. As a case study we consider lattices of coupled particles with nonlinear nearest neighbour interaction and prove the existence of exponentially spatially localised breathers exhibiting either even-parity or odd-parity symmetry under necessary non-resonant conditions accompanied with the proof of energy bounds of the solutions. 
\end{abstract}

\maketitle

\noindent 
\section{Introduction}
The study of nonlinear lattice dynamical systems is an active area of research of continuously growing interest.  For the motivation, efforts, and research achievements in this major field of nonlinear science we refer to \cite{BK}-\cite{TD2}. A turning point can be traced back to  the end of 30’s of the previous century when the striking studies of Frenkel and Kontorova \cite{FK1} appeared.  These studies on crystal dislocations brought about the derivation of the Discrete Klein-Gordon equation (DKG). In the 50’s came a revolution with the famous Fermi-Pasta-Ulam-Tsingou (FPUT) experiment \cite{FPUT1}-\cite{TD2}, investigating the equipartition of energy in a weakly nonlinear system of coupled oscillators. The unexpected and surprising results of this experiment are identified with the birth of the mathematics and physics of solitons. The FPUT system, the DKG equation, the Discrete Nonlinear Schrödinger equation (DNLS) and their variants are the prototype systems for the study of localised phenomena (in crystals, nonlinear-waveguide arrays, Bose-Einstein condensates, DNA macromolecules-Davydov mechanism).

Due to their relevance to the above mentioned nonlinear phenomena of restrained energy redistribution, the problem of the existence of discrete localised modes in the form of solitons and breathers in the above lattice systems has attracted tremendous interest,  not only from the physical but also from the  rigorous mathematical viewpoint.  Results regarding  proofs of existence and stability of discrete solitons and breathers are based on a variety of mathematical methods: on the one hand, integrability theory applies to only very few  systems such as the Toda and the Ablowitz-Ladik lattices \cite{Toda1}-\cite{AL2}, while the wide class of non-integrable systems has been studied with  the  continuation of (trivially localised) solutions from the anti-integrable limit
for sufficiently weak coupling strengths,  Lyapunov-Schmidt decomposition, normal form and center manifold reduction methods, Floquet monodromy theory, as well as, direct and minimax variational approaches. We refer to \cite{MacKay}-\cite{APB} and cited works therein.

In this work, pursuing further our studies initiated in \cite{JFPT} for DKG-chains, we claim that a fixed point theory approach, based on the celebrated Schauder's fixed point theorem \cite{Schauder},  can offer  a comparatively concise and to our knowledge, novel approach  compared to the aforementioned methods, for the proof of existence of localised waveforms for nonlinear lattices. In the present paper, we facilitate this fixed point approach to a general
class of  lattices with nonlinear interaction between its constituents 
and with general potentials whose second derivative is bounded by  a power-law growth. Prominent examples are the $\alpha$-FPUT and $\beta$-FPUT systems or combinations of both. 
Other important examples of interaction potentials include  periodic and  bounded soft potentials. We think that one of the main assets of the Schauder's fixed point method is that it  establishes, apart from existence, directly a requested degree of spatial localisation \cite{Pankov2}, \cite{Zhang3}, \cite{Pankov3}. That is, the breathers are proved to be exponentially localised as it is the case in  physically significant applications  (see \cite{PRL1}, \cite{PhysRep1} and references therein). Another potential advantage, is that the method is not restricted to the existence of small amplitude oscillations \cite{James}, and is promising for its  extension to higher dimensional set-ups (where localisation effects are more intricate \cite{FG}), as its major characteristics (particularly the derivation of suitable estimates), seem not to be affected by higher-dimensionality. 

The presentation of the paper is as follows. In section \ref{sec2} we discuss the model systems and main assumptions on the interaction potential.  Section \ref{sec3} contains the description and properties of the functional setting essentially involving weighted sequence spaces and the proof of the main result. We claim that the whole approach can be of general  interest: Firstly, it implements a Fourier series representation of the solutions to establish compactness embedding properties of the considered functional spaces; these can be useful generically for the analysis of lattice dynamical systems (existence and dynamics, see \cite{NT2005},\cite{DNat} and references therein). Secondly,  facilitating the Schauder's fixed point approach and investigating the contraction regime for the nonlinear map in this set-up, gives rise to physically relevant and consistent restrictions on the characteristics of the breather solutions, namely, their frequency (natural non-resonant conditions) and their magnitude measured 
 in the suitable weighted norm (norm  upper and lower bounds); importantly, the latter exhibits dependence on the localisation degree (as dictated by the relevant weighted norm) and the growth of the nonlinearity for the general class of potentials considered herein (see Theorem \ref{mainT} and subsection \ref{RemF}).
Finally, section \ref{sec4} summarises the main findings and gives an outlook to  extensions for future studies.  

\section{Description of the system}
\label{sec2}
The  system of interacting particles on 
infinite  lattices is given by the following set of coupled equations
\begin{equation}
\frac{d^2 q_n}{dt^2}=-V^{\prime}(q_n-q_{n-1})+
V^{\prime}(q_{n+1}-q_n),\,\,\, n\in {\mathbb{Z}},\label{eq:system}
\end{equation}
where $q_n$ is the displacement of a particle from its equilibrium position.  
The function $V(x)$ is the interaction potential between nearest neighbours  and 
the prime $^{\prime}$  stands for the derivative with respect to the argument $x$.
This system has a Hamiltonian structure related to the energy
\begin{equation}
 H=\sum_{n\in {\mathbb{Z}}}\left(\frac{1}{2}p_n^2+V(q_n-q_{n-1})\right),\nonumber
\end{equation}
 and it is time-reversible with respect to the involution $p\mapsto -p$. 
 In addition to the energy $H=E$, the total momentum $P=\sum_n p_n$ is conserved too.
 
 Discrete breathers can be characterised as follows:
 \begin{eqnarray}
  q_n(t+T)&=&q_n(t),\,\,\,p_n(t+T)=p_n(t),\,\, n\in {\mathbb{Z}},\label{eq:periodic}\\
  \lim_{n\rightarrow \pm \infty}q_n&=&c_{\pm},\,\,\, \lim_{n\rightarrow \pm \infty}p_n=0,\label{eq:limits}
 \end{eqnarray}
 and the difference $c_+ -c_{-}$ is the lattice deformation induced by the breather. 
 Furthermore,  breather solutions are characterised by zero time average, i.e.
 \begin{equation}
  \int_0^T q_n(t)dt=0,\,\,\,n\in {\mathbb{Z}}.\label{eq:zeroaverage}
 \end{equation}

We consider periodic solutions exhibiting even-parity and odd-parity symmetry.
To this end we introduce the following linear mappings:
\begin{eqnarray}
 R_{odd}\,:\,(b)_{n \in \mathbb{Z}} \mapsto (b)_{n \in \mathbb{Z}},\,\,\,\left(R_{odd}\cdot b\right)_n&=&-b_{-n},\,\,\,b_{n \in \mathbb{Z}}\in \mathbb{R},\nonumber\\
 R_{even}\,:\,(b)_{n \in \mathbb{Z}} \mapsto (b)_{n \in \mathbb{Z}},\,\,\,\left(R_{even}\cdot b\right)_n&=&-b_{-(n+1)},\,\,\,b_{n \in \mathbb{Z}}\in \mathbb{R}.\nonumber
\end{eqnarray}
A periodic solution $q(t)=q(t+T)$ has  even-parity symmetry if
\begin{equation}
 R_{even}\cdot q(t)=q(t),\,\,\,\forall t \in \mathbb{R},\label{eq:Reven}
\end{equation}
and possesses odd-parity symmetry if
\begin{equation}
 R_{odd}\cdot q(t+T/2)=q(t),\,\,\,\forall t \in \mathbb{R}.
\end{equation}
The even-parity symmetric solution is characterised by 
 amplitude profiles centered between  sites $n=-1$ and $n=0$ while the odd-parity symmetry solution belongs to amplitude profiles centered at site $n=0$. 
 
\vspace*{0.5cm}

The existence theorem will be proved under the following assumption:
\begin{itemize}
	\item
\textbf{A:}  
The interaction potential $V(x)$ has a global minimum at $x=0$ and   is of  such form that
\begin{equation}
 V(x)=\frac{1}{2}x^2+W(x),\label{eq:Af}
\end{equation}
with $W:\,\,{\mathbb{R}}\,\rightarrow\,{\mathbb{R}}$ at least twice continuously 
differentiable and  $W(0)=W^\prime(0)=0$. $W(x)$ is of the form that $V^{\prime \prime}(0)=a>0$ holds.
Further, we assume that $W$ satisfies for some constant $\overline{K}>0$ and $\alpha\geq 0$, 
such that
\begin{equation}
|W''(x)|\leq \overline{K}|x|^\alpha,\;\;\forall x\in\mathbb{R}.
	\label{eq:assumptions1}
\end{equation}
%
\end{itemize}

 Examples of interaction potentials which satisfy condition $\mathbf{A}$ are given by the $\alpha$-FPUT and $\beta$-FPUT system with  $W(x)=\alpha x^3/3$ and $W(x)=\beta x^4/4$ as well as combinations of both
\cite{Wattis},\cite{Arioli}.

The solutions of the system obtained when linearising  equations (\ref{eq:system}) around the equilibrium $q_n=0$, are superpositions of plane wave solutions (phonons) 
\begin{equation}
 q_n(t)= \exp(i(kn-\omega t)),\nonumber
\end{equation}
with frequencies 
\begin{equation}
 \omega^2(k) =4\sin^2\left(\frac{k}{2}\right),\,\,\,k\in[-\pi,\pi].\nonumber
\end{equation}
These (extended) states disperse. Therefore, the  frequency $\Omega$ of a  localised time-periodic solution must satisfy the non-resonance condition $\Omega \neq |\omega(k)|/m$ for any integer $m\ge 1$. This requires 
$\Omega^2 >4$ as a necessary condition for the existence of localised time-periodic solutions of system (\ref{eq:system}).

For the forthcoming analysis, it is convenient to pass to relative variables 
defined by
\begin{eqnarray}
 x_n&=&q_n-q_{n-1},\label{eq:newx}\\
 p_n&=&y_n-y_{n+1},\label{eq:newp}
\end{eqnarray}
where the original variables satisfy
\begin{equation}
 \sum_{n\in \mathbb{Z}}|p_n|<\infty,\,\,\,\sum_{n\in \mathbb{Z}}|q_n-q_{n-1}|<\infty.
\end{equation}
which imply that
\begin{equation}
 \sum_{n\in \mathbb{Z}}|y_n-y_{n+1}|<\infty,\,\,\, \sum_{n\in \mathbb{Z}} |x_n|<\infty.
\end{equation}
From (\ref{eq:limits}), it follows that $\lim_{n\rightarrow -\infty}y_n=0$ so that the change of variables  (\ref{eq:newx})-(\ref{eq:newp}) is invertible. That is, for given $\{(x_n,y_n)\}_{n\in \mathbb{Z}}$, one gets the relations
\begin{equation}
 p_n=y_n-y_{n+1},\,\,\,q_n=\sum_{l=-\infty}^n x_l.
\end{equation}
The transformation (\ref{eq:newx})-(\ref{eq:newp}) is symplectic, because
\begin{equation}
 \sum_n dp_n\wedge dq_n=\sum_n\left(dy_n\wedge dq_{n-1}+dy_n \wedge dx_n-dy_{n+1}\wedge dq_n\right)=\sum_ndy_n\wedge dx_n,\label{eq:symplectic}
\end{equation}
and from the Hamiltonian
\begin{equation}
 H=\sum_n \left(\frac{1}{2}(y_n-y_{n+1})^2+V(x_n)\right),
\end{equation}
one derives the canonical equations 
\begin{equation}
 \dot{x}_n=2y_n-y_{n+1}-y_{n-1},\,\,\,\dot{y}_n=-V^\prime(x_n),
\end{equation}
or equivalently
\begin{equation}
 \ddot{x}_n=V^\prime(x_{n+1})+V^\prime(x_{n-1})-2V^\prime(x_{n}).\label{eq:transformed}
\end{equation}
In (\ref{eq:symplectic}) the symbol $\wedge$ denotes the  standard wedge product.
The relations (\ref{eq:periodic})-(\ref{eq:zeroaverage}), satisfied by breather solutions of the original system (\ref{eq:system}), 
when passing to the new variables, read as
 \begin{eqnarray}
  x_n(t+T)&=&x_n(t),\,\,\,y_n(t+T)=y_n(t),\,\, n\in {\mathbb{Z}}\nonumber\\
  \lim_{n\rightarrow \pm \infty}x_n&=&0,\,\,\,\lim_{n\rightarrow \pm \infty}y_n=0,\nonumber\\
    \int_0^T x_n(t)dt&=&0,\,\,\,n\in {\mathbb{Z}}.\nonumber
\end{eqnarray}

\section{Existence of exponentially localised breathers}
\label{sec3}
\setcounter{equation}{0}
In the following we prove 
the existence  of  localised periodic 
solutions of system (\ref{eq:system}) on the infinite one-dimensional lattice.
To this end, some appropriate function spaces are introduced on which the original 
problem is presented as a fixed point problem for a corresponding operator. 
Utilising Schauder's 
Fixed Point Theorem, we establish the existence of exponentially localised solutions, extending the approach followed in \cite{JFPT}, to a suitable functional spaces- setup.

To establish the required spatial localisation of the solutions, we introduce suitable weighted function spaces.
First we consider the exponentially weighted Hilbert space of square-summable sequences, ${l}^2_w(\mathbb{Z})$,  defined as
\begin{equation}
	l^2_w=\bigg\{u_n\in \mathbb{R}\,:\,||u||_{l_w^2}^2:=
	\sum_{n}\,w_n|u_n|^2\,\bigg\}\,.\label{eq:norm}
\end{equation}
and we use 
$w_{n \in \mathbb{Z}}=\exp(\lambda |n|)$ and $\lambda \ge 0$ complying according to (\ref{eq:Reven})  and (\ref{eq:newx}) with the even-parity amplitude profile.
  Denoting by $u^{even}$ and $u^{odd}$ the even and odd parity amplitude pattern, respectively, 
  it can be readily shown that 
  \begin{equation}
   ||u^{odd}||_{l^2_{w}}<||u^{even}||_{l^2_{w}}.\label{eq:normodd}
  \end{equation}

Then,  we denote by 
\begin{equation}
	{X}_{0} =  \bigg\{ u\in L^{2}_{per}([0,T];{l}_w^2)\,:\,\int_0^T u_n(t)dt=0,\,\,\, n\in {\mathbb{Z}}
	\bigg\},\nonumber
\end{equation}
the space of $T-$periodic square integrable
functions in time with zero time average, with values in $l^2_w$. Evidently, $X_0$ is a closed convex subspace of $L^{2}_{per}([0,T];{l}_w^2)$. We also consider the Sobolev space 
\begin{equation*}
	X_2=\bigg\{u\in H^{2}_{per}([0,T];{l}_w^2)\bigg\},
\end{equation*}
containing  $T-$periodic 
functions of time  assuming values in $l_w^2$ which, 
together with their weak derivatives  up to second order are in $X_0$. The above spaces are endowed with the following norms:
\begin{eqnarray*}
	|| u ||_{X_0}^2&=&\frac{1}{T}\int_{0}^{T}\,|| u(t)||^2_{{l}_w^2}dt,
	\nonumber\\
	|| u ||_{X_2}^2&=&\frac{1}{T}\int_{0}^{T}\,\left(|| u(t)||^2_{{l}_w^2}
	+|| D{u}(t)||^2_{{l_w}^2}+|| D^2{u}(t)||^2_{{l^2_w}}\right)dt.\nonumber
\end{eqnarray*} 
For an element $u\in X_2$, we may consider the  Fourier-series expansion of $u_n(t)$ with respect to time $t$ and space variable $n$,  determined by 
\begin{eqnarray}
	\label{eq:FTX1}
	u_n(t)&=&\sum_{m\in {\mathbb{Z}}\setminus \{ 0 \}}
	\hat{u}_{n,m} \exp\left( i\Omega m t\right),\,t\in [0,T],\,\,\,\,\,\,{\hat{u}}_{n,m}=\frac{1}{T} \int_{0}^{T} u_n(t)\exp(-i \Omega m t)dt,\\
	\hat{u}_{n,m}&=&\overline{\hat{u}}_{n,-m},
	\label{eq:FTX2}
\end{eqnarray}
and 
\begin{eqnarray}
	\label{eq:FTspatial0}
	u_n(t)&=&\frac{1}{2\pi}\int_{0}^{2\pi}\tilde{u}_k(t)\exp(ikn)dk,\,\,\,\,\,\tilde{u}_k(t)=\sum_{n}u_n(t)\exp(-i\,k\,n),\\
\label{eq:FTspatial}	
	\tilde{u}_{k}&=&\overline{\tilde{u}}_{-k}.
\end{eqnarray}
Then, using \eqref{eq:FTX1}-\eqref{eq:FTspatial}, we may have the representation of $u$, given by
\begin{equation}
	u_n(t)=\sum_{m\in {\mathbb{Z}}\setminus \{0\}}\frac{1}{2\pi} \int_{0}^{2\pi} {\hat{\tilde{u}}}_{k,m}\exp(ikn)dk\exp(i\Omega mt).\label{eq:FTst}
\end{equation}
The following compactness lemma will be particularly useful for our approach.
\begin{lemma}
\label{comlemma}
$X_2$ is compactly embedded in $X_0$ (denoted as $X_2 \Subset X_{0}$).
\end{lemma}
\textbf{Proof:} The continuous embedding $X_2\subseteq X_0$ is clear due to the inequality
\begin{equation}
	\label{conembe}
	||u||_{X_0}\le ||u||_{X_2}.\nonumber
\end{equation}
Now, in order to show the compactness of the embedding, let $(u_n)_{n \in \mathbb{N}}\in X_2$ be a bounded sequence, i.e. $||u_n||_{X_2}\le M<\infty$ for all $n \in \mathbb{N}$. Then, there exists a weakly convergent subsequence (not relabeled) $u_n\rightharpoonup u\in X_2$. The continuous embedding  $X_2\subseteq X_0$ implies 
the weak convergence $u_n\rightharpoonup u$
in $X_0$. 
For the strong convergence $u_n\rightarrow u$ in $X_0$,
considering the difference $||u_n-u||_{X_0}^2$, we have the estimate:
\begin{eqnarray}
	||u_n-u||_{X_0}^2&=&\frac{1}{T}\int_0^T||u_n(t)-u(t)||_{l^2_w}^2 dt\nonumber\\
	&=& \frac{1}{T}\int_0^T\sum_{k\in \mathbb{Z}}w_k| u_n^k(t)-u^k(t)|^2 dt\nonumber\\
	&=& \frac{1}{T}\int_0^T\sum_{|k|\le K}w_k| u_n^k(t)-u^k(t)|^2 dt+ \frac{1}{T}\int_0^T\sum_{|k|>K} w_k|u_n^k(t)-u^k(t)|^2 dt\nonumber\\
	&=&\frac{1}{T}\int_0^T\sum_{|k|\le K}w_k| u_n^k(t)-u^k(t)|^2 dt+\sum_{|k|>K}w_k\sum_{m\neq 0}|\hat{u}_{n,m}^{k}-\hat{u}_m^{k}|^2\nonumber\\
	&\le& \frac{1}{T}\int_0^T\sum_{|k|\le K}w_k| u_n^k(t)-u^k(t)|^2 dt+\sum_{|k|>K}w_k\sum_{m\neq 0}|\hat{u}_{n,m}^{k}|^2+\sum_{|k|>K}w_k\sum_{m\neq 0}|\hat{u}_m^{k}|^2\nonumber\\
	&\le& \frac{1}{T}\int_0^T\sum_{|k|\le K}w_k| u_n^k(t)-u^k(t)|^2 dt+\sum_{|k|>K}w_k|{u}_n^{k}|^2\sum_{m\neq 0}\frac{1}{m^4}+\sum_{|k|>K}w_k|{u}^{k}|^2\sum_{m\neq 0}\frac{1}{m^4}\nonumber\\
	&=& \frac{1}{T}\int_0^T\sum_{|k|\le K}w_k| u_n^k(t)-u^k(t)|^2 dt+2\zeta (4)\left( \sum_{|k|>K}w_k|{u}_n^{k}|^2+\sum_{|k|>K}w_k|{u}^{k}|^2 \right),\label{eq:embedding1}
\end{eqnarray}
where $\zeta(s)=\sum_{m=1}^{\infty}m^{-s}$ is the Riemann zeta function, $({u}_n^{k})_{k\in I}\in l^2_w$,  
$({u}^{k})_{k\in I}\in l^2_w$, $I=(-K,...,0,...K)$.  We  also used that since $u\in H^2_{per}[0,T]\Rightarrow |\hat{u}_m|\le C/{|m|^2}$.
Then for sufficiently large $K$ we get  
\begin{eqnarray}
	2\zeta (4)\left( \sum_{|k|>K}w_k|{v}_n^{k}|^2+\sum_{|k|>K}w_k|{v}^{k}|^2 \right)
	&\le& \frac{2\epsilon}{3}.\label{eq:zeta}
\end{eqnarray}

Exploiting that weak and strong convergence coincide in  finite-dimensional spaces, by choosing $n$ sufficiently large, we get for the first term on the right-hand side of
(\ref{eq:embedding1}), that 
\begin{equation}
	\frac{1}{T}\int_0^T\sum_{|k|\le K}w_k| u_n^k(t)-u^k(t)|^2 dt\le \frac{\epsilon}{3}.\label{eq:weak} 
\end{equation}
Combining (\ref{eq:zeta}) and (\ref{eq:weak}) we get
\begin{equation}
	||u_n-u||_{X_0}^2\le \epsilon,
\end{equation}
and the proof is finished.

\  \ $\Box$

Let us also recall the two versions of the Schauder's fixed point Theorem \cite{Schauder}.
\begin{theorem}
\label{SFPT}
\begin{enumerate}
	\item (First version): Let $G$ be a closed bounded convex subset of a Banach space $X$.  Assume that $f:\,G\mapsto G$ is compact. Then, $f$ has at least one fixed point in $G$.
\item (Second version):  Let $G$ be a closed convex subset of a
	Banach space $X$ and $f$ a  continuous map of $G$ into a compact subset of $G$. Then,  $f$ has at least one  fixed point.
\end{enumerate}
\end{theorem}

With the preparations above, we may proceed to the statement and proof of the main result.

\begin{theorem} 
\label{mainT}	(i). Let condition {\bf A}  hold and suppose 
\begin{equation}
\Omega^2> 4. \label{eq:nonresonance0}
\end{equation}
Then 
there exist nonzero sequences  $x^{even}\equiv \{x_n^{even}\}_{n\in {\mathbb{Z}}} 
\in H^2([0,T];l^2)$ and $x^{odd}\equiv \{x_n^{odd}\}_{n\in {\mathbb{Z}}} 
\in H^2([0,T];l^2)$ with zero time average $\int_0^T x_n(t)dt=0$ and $||x||_{X_0}\le R$, 
where 
\begin{equation}
R\le \left(\frac{\Omega^2-4}{\overline{K}\,\sqrt{2(1+\cosh(\lambda))}}\right)^{1/\alpha}:=R_{\mathrm{max}}.  \label{eq:As0}
\end{equation}	
These solutions are exponentially localised time-periodic solutions of system (\ref{eq:transformed}) with period   $T=\frac{2\pi}{\Omega}$ possessing even-parity and odd-parity symmetry, respectively. \\
(ii). Assume further that
	\begin{equation}
	\Omega^2>4+\overline{K}\,\sqrt{2(1+\cosh(\lambda))}.\label{eq:nonresonance}
	\end{equation}
Then, the corresponding solutions satisfy the lower and upper bounds in their $X_0$-norm	
	\begin{equation}
	R_{\mathrm{crit}}:=\left(\frac{\Omega^2-4}{\overline{K}\,\sqrt{2(1+\cosh(\lambda))}}\right)^{\frac{1}{1+\alpha}}\le R\le R_{\mathrm{max}}.  \label{eq:As}
	\end{equation}
\end{theorem}
{\bf Proof:} For part (i), we prove the existence of solutions when $\Omega^2>4$, which is a natural condition for the existence of localised modes. We will show that the solutions are satisfying the upper bound  \eqref{eq:As0},  by implementing  the two versions of the Schauder's fixed point Theorem \ref{SFPT} to the suitably defined map. For part (ii), assuming \eqref{eq:nonresonance}, we examine the contraction mapping regime for this map, identifying a ring in $X_0$ of nontrivial solutions satisfying the lower bound of \eqref{eq:As}.\\
\textit{(i): Existence by the Schauder's fixed point Theorem of even-parity and odd-parity solutions when \eqref{eq:nonresonance0} holds}. 
First we show the existence of even-parity symmetry breathing solutions. Afterwards we argue that  by virtue of relation (\ref{eq:normodd}) the existence proof applies to the odd-parity breathing mode too.
 For the ease of notation we denote the solution by $x$ instead of $x^{even}$. 
We shall provide two alternatives of the proof by implementing  the two versions of the Schauder's fixed point Theorem \ref{SFPT}.  For this purpose, it is  convenient to rewrite Eq.\,(\ref{eq:transformed})  as:
\begin{equation}
 \ddot{x}_n-[x_{n+1}-2x_n+x_{n-1}]=W^{\prime}(x_{n+1})+W^{\prime}(x_{n-1})-2W^{\prime}(x_{n}).
\label{eq:ref1}
\end{equation}
Thus, only the right-hand side  of  (\ref{eq:ref1}) features terms nonlinear in $x$. We shall express  
(\ref{eq:ref1}) as a fixed point equation in $x$.\\
{\em I. $1^{st}$ version of the proof:} We relate the left-hand side of (\ref{eq:ref1}) to  the 
linear mapping: $M\,:\,X_2\,\rightarrow\,X_{0}$:
\begin{equation}
 M(u_n)=\ddot{u}_n- [u_{n+1}-2u_n+u_{n-1}]. \nonumber 
\end{equation}
Then, applying the operator $M$ to the Fourier elements $\exp(ikn)\exp(i\,\Omega m t)$ 
in the representation (\ref{eq:FTst}), we get that
\begin{equation}
 M\exp(ikn)\exp(i\,\Omega mt)=\nu_m(k) \exp(ikn)\exp(i\,\Omega mt),\nonumber
\end{equation}
where 
\begin{equation}
 \nu_m(k)= - \Omega^2\,m^2+4\sin^2\left(\frac{k}{2}\right),\,\, m\in \mathbb{Z} \setminus \{ 0 \},\,\,k\in [0,2\pi].\nonumber
\end{equation}
As by assumption (\ref{eq:nonresonance}) $\Omega^2>4$, it is guaranteed that $\nu_m(k) \ne 0$, for all $m\in {\mathbb{Z}} \setminus \{ 0 \}$ and for all $k\in [0,2\pi] $, the mapping $M$ possesses an inverse $M^{-1}$, obeying 
$M^{-1}\exp(i\,(\Omega mt+kn))=(1/\nu_m(k))\exp(i\,(\Omega mt+kn))$.  
For the norm of the linear operator 
$M^{-1}:\,\,X_{0} \rightarrow X_2$ 
one derives  
the upper bound:
  \begin{eqnarray}
  || M^{-1} ||_{X_{0},X_2}^2&=&
   \sup_{_{||u||_{X_0}=1}}|| M^{-1}\,u ||_{X_2}^2\nonumber\\
   & =&\sup_{_{||u||_{X_0}=1}}\frac{1}{T}\int_0^T\,\left[|| M^{-1} u(t)||^2_{{l}_w^2}
   +|| D{M^{-1}u(t)} 
   ||^2_{{l}_w^2}+|| 
   D^2{M^{-1}u(t)}
   ||^2_{{l}_w^2}\right]dt\nonumber\\
  &=&\sup_{||u||_{X_0}=1}\frac{1}{T}\int_0^T\,\sum_{n\in \mathbb{Z}}
  w_n \left(\left| \sum_{m^\prime}
\frac{1}{2\pi} \int_{0}^{2\pi}\frac{\hat{\tilde{u}}_{k,m}}{\nu_m(k)}\exp(ikn)\exp(i\Omega m t) dk \right|^2\right.\nonumber\\   
  &+&\left|\sum_{m^\prime}\frac{1}{2\pi} \int_{0}^{2\pi}
    \frac{i \Omega m \hat{\tilde{u}}_{k,m}}{\nu_m(k)} \exp(ikn)\exp(i\Omega m t) dk \right|^2 \nonumber\\
 &+&\left.\left|\sum_{m^\prime}\frac{1}{2\pi} \int_{0}^{2\pi}
    \frac{(i \Omega m)^2 \hat{\tilde{u}}_{k,m}}{\nu_m(k)} \exp(ikn)\exp(i\Omega m t) dk \right|^2 \right)dt\nonumber\\   
&\le&\sup_{m\in {\mathbb{Z}}\setminus \{ 0 \}} \sup_{k\in [0,2\pi]} \frac{1
    +(\Omega m)^2+(\Omega m)^4}{|\nu_m(k)|^2} \nonumber\\
         &\cdot & \sup_{||u||_{X_0}=1} \frac{1}{T} \int_0^T\,\sum_{n\in \mathbb{Z}}
    w_n \left| \sum_{m^\prime}
  \frac{1}{2\pi} \int_{0}^{2\pi} \hat{\tilde{u}}_{k,m} \exp(ikn)\exp(i\Omega m t) dk \right|^2 dt\nonumber\\
&\le&
\frac{1+\Omega^2+\Omega^4}{(\Omega^2-4)^2} \sup_{||u||_{X_0}=1} ||u||_{X_0}^2
  =  \frac{1+\Omega^2+\Omega^4}{(\Omega^2-4)^2}\le \frac{(1+\Omega^2)^2}{(\Omega^2-4)^2}<\infty,
   \label{eq:boundL}
  \end{eqnarray}
verifying the boundedness of $M^{-1}$   
and we used the notation $\sum_{l^\prime} = \sum_{l\in {\mathbb{Z}}\setminus \{ 0 \}}$. 
When we consider $M^{-1}$ as a map $M^{-1}:\,\,X_{0} \rightarrow  X_0$,  we note that 
 \begin{equation}
  || M^{-1} ||_{{X}_{0},{X}_0}\le \frac{1}{\Omega^2-4}.
  \label{eq:boundL0}
 \end{equation}
 
In order to deal with the nonlinear terms, we attribute to the  right-hand side  of (\ref{eq:ref1}) 
the nonlinear operator $N\,:\,X_{0}\,\rightarrow\,X_{0}$, as
\begin{equation}
 (N(u))_n= W^{\prime}(u_{n+1})+W^{\prime}(u_{n-1})-2W^{\prime}(u_n).\nonumber
\end{equation}
To establish  that the  operator $N$ is continuous on $X_{0}$, we prove that  
 $N$ is Frechet differentiable at any $u$, with bounded derivative. We observe that
\begin{equation}
 N^{\prime}(u):\,h\in X_{0}\,\mapsto N^{\prime}(u)[h] =W^{\prime\prime}(u)h \in X_{0},\nonumber
\end{equation}
where $(W^{\prime\prime}(u))_n=W^{\prime\prime}(u_{n+1})+W^{\prime\prime}(u_{n-1})-2W^{\prime\prime}(u_n)$. Therefore, we derive the estimate
\begin{eqnarray}
 || N^{\prime}(u)[h]||_{X_{0}}^2&=& \frac{1}{T}
 \int_0^T
 \sum_{n \in \mathbb{Z}}{w_n} \left|\,\left(W^{\prime\prime}(u_{n+1}(t))+W^{\prime\prime}(u_{n-1}(t))-2W^{\prime\prime}(u_{n}(t))\right)\,h_n(t)\right|^2\,dt\nonumber\\
   &\le& 
  \frac{1}{T} \int_0^T
  \sum_{n \in \mathbb{Z}} w_n \overline{K}^2(|u_{n+1}(t)|^{2\alpha}+|u_{n-1}(t)|^{2\alpha}+2|u_{n}(t)|^{2\alpha})\,|h_n(t)|^2\,dt\nonumber\\
&\le& 4\overline{K}^2 \sup_{n\in {\mathbb{Z}}} \max_{t\in [0,T]} |u_n(t)|^{2\alpha}\,|| h||^2_{X_0}=A^2 || h||^2_{X_0},\nonumber
\end{eqnarray}
where $A^2=4\overline{K}^2 \sup_{n\in {\mathbb{Z}}} \max_{t\in [0,T]} |u_n(t)|^{2\alpha} $. That is
\begin{equation}
 || N^{\prime}(u)||_{{\cal{L}}(X_{0},X_{0})}\le A,\label{eq:Frechet}
\end{equation}
implying the (uniform) boundedness of the differential. Let us now use  as the closed convex subset $Y_0$ of $X_0$, its closed ball centered at $0$ of radius $R$,
\begin{equation}
 {Y}_{0} =  \left\{ u\in X_0\;\;:\;\;|| u||_{X_0}\le R
 \right\}.\nonumber
\end{equation}
For even-parity solutions we get for the range of $N$ on $Y_0$ the bound
\begin{eqnarray}
 || N(u) ||_{X_{0}}^2&=&\int_0^T
 \sum_{n \in \mathbb{Z}}{w_n}|   W^{\prime}(u_{n+1}(t))+W^{\prime}(u_{n-1}(t))-2W^{\prime}(u_{n}(t))        |^2\,dt\nonumber\\
 &\le& 
 \int_0^T
 \sum_{n \in \mathbb{Z}}{w_n} \left(|W^{\prime}(u_{n+1}(t))|^2+|W^{\prime}(u_{n-1}(t))|^2+2|W^{\prime}(u_n(t))|^2\right)\,dt\nonumber\\
 &\le&
 \overline{K}^2\int_0^T\sum_{n \in \mathbb{Z}}\left(\frac{w_{n}}{w_{n+1}}w_{n+1}|u_{n+1}(t)|^{2(\alpha+1)}+\frac{w_{n}}{w_{n-1}}w_{n-1}|u_{n-1}(t)|^{2(\alpha+1)}+2w_n|u_n(t)|^{2(\alpha+1)}\right)\,dt\nonumber\\
 &\le& 2\overline{K}^2(1+\cosh(\lambda))||u||_{\tilde{X}_0}^{2(\alpha+1)}\le 2\overline{K}^2(1+\cosh(\lambda))||u||_{{X}_0}^{2(\alpha+1)}\nonumber\\
 &\le & 2\overline{K}^2  (1+\cosh(\lambda)) R^{2(\alpha+1)},\qquad \forall u\in Y_0,
 \label{eq:rangeN}
\end{eqnarray}
where $\tilde{X}_0=\{u\in L_{per}^2([0,T];l_w^{2(\alpha+1)})\,:\,\int_0^T u_n(t)dt=0,\,\, n\in {\mathbb{Z}}\}$ and we used $X_0\subset \tilde{X}_0$.
(A similar expression is obtained for odd-parity solutions.)

In order to apply Theorem \ref{SFPT}, the final step is to express the problem  (\ref{eq:ref1})
as a fixed point equation in terms of a mapping $Y_0\,\rightarrow\, Y_{0}$:
\begin{equation}
 x=M^{-1}\circ N(x)\equiv S(x).\label{eq:compose}
\end{equation}
Clearly, $S$ is continuous on $X_{0}$, as the relations (\ref{eq:boundL}) and (\ref{eq:Frechet}) establish that its constituents $M^{-1}$ and $N$ are continuous.
Next, we show that $S({Y}_0)\subseteq {Y}_0$.  Using (\ref{eq:boundL}) and (\ref{eq:rangeN}), we have
\begin{eqnarray}
\label{Ring1}
 || S(x)||_{{X}_0}&=&|| M^{-1}(\,N(x))||_{{X^0}}
 \le || M^{-1} ||_{{X}_0,X_0}\cdot||  \,N(x)||_{{X}_0}\nonumber\\
 & \le& 
 \frac{ \sqrt{2(1+\cosh(\lambda)}\,\, \overline{K}  
}{\Omega^2-4} R^{\alpha+1}\le R,\,\,\, \forall x\in Y^0,
\end{eqnarray}
assuring by assumption (\ref{eq:As}), that indeed
\begin{equation}
S({Y}_0)\subseteq {Y}_0.\nonumber
\end{equation} 
As $M^{-1}$ maps $X_0$ to $X_2 \Subset X_{0}$, one has 
$S(Y_0) \subseteq Y_0 \cap X_2$. Therefore the map $S=M^{-1}\circ N$, viewed as a map $S: Y_0\subseteq X_0\mapsto Y_0 \subseteq X_0$, is compact. The first version of Schauder’s fixed point theorem
implies then that the fixed point equation $x =  S(x)$ has at least one solution of even-parity symmetry. 
 In addition, due to the relation (\ref{eq:normodd}), the above existence proof is valid for the odd-parity symmetry solution when $\Omega$ is considered as its associated breathing frequency.
\\
{\em II. $2^{nd}$ version of the proof:} For the application of the second version of Schauder's fixed point theorem, we verify that the range of $S$ is contained in a compact subset of $Y_0$. 
Representing  $N(u) \in Y_0$ in terms of its spatial and temporal Fourier-transforms  as
\begin{equation}
	(N(u))_n(t)=\sum_{m\in {\mathbb{Z}}\setminus \{ 0 \}} \frac{1}{2\pi} \int_{0}^{2\pi}
	\hat{\tilde{N}}_{k,m}\exp(ikn)dk \exp\left(i\Omega m t\right),\label{eq:FT}
\end{equation}
the Fourier coefficients of $M^{-1}(N(u))$ fulfill for all $u\in Y_0$, for all $k\in[0,2\pi]$, and for all $m\in {\mathbb{Z}}\setminus \{ 0 \}$
\begin{eqnarray}
	& &\left|\frac{\hat{{\tilde{N}}}_{k,m}}{- \Omega^2\,m^2+4 \sin^2\left(\frac{k}{2}\right)}\right|^2\le 
	\frac{\sum_{m\in {\mathbb{Z}}\setminus \{ 0 \}}
		\frac{1}{2\pi} \int_{0}^{2\pi}\left|\hat{{\tilde{N}}}_{k,m}\right|^2dk}{m^4(\Omega^2-4)^2}\nonumber\\
	&=&\frac{\parallel N(u) \parallel_{{X}_{0,0}}^2}{m^4(\Omega^2-4)^2},\label{eq:FcN}
\end{eqnarray}
where $X_{0,0}$ is determined by 
\begin{equation}
	{X}_{0,0} =  \left\{ u\in L^{2}_{per}([0,T];{l}^2)\,\vert\,\int_0^T u_n(t)dt=0,\,\,\, n\in {\mathbb{Z}}
	\right\}.\nonumber
\end{equation}
Since $l_w^2\subset l^2$ and the inequality $||u||_{l^2} \le ||u ||_{l_w^2}$ holds,  we get
\begin{eqnarray}
	& &\left|\frac{\hat{{\tilde{N}}}_{k,m}}{- \Omega^2\,m^2+4 \sin^2\left(\frac{k}{2}\right)}\right|^2\le \frac{2{\overline{K}}^2(1+\cosh(\lambda))R^{2(\alpha+1)}}{m^4(\Omega^2-4)^2}\le \left(\frac{R}{m^2}\right)^2.
\end{eqnarray} 
Hence, we conclude that $S$ maps $Y_0$ into the subset 
\begin{equation}
	Y_{0,c}=\left\{u=\left\{u\right\}_{n\in {\mathbb{Z}}} \in Y_0,\,u_n(t)=\sum_{m\in {\mathbb{Z}}\setminus \{ 0 \}}
	\frac{1}{2\pi} \int_{0}^{2\pi}\hat{{\tilde{u}}}_{k,m}\,\exp(ikn)dk\exp\left( i \Omega m t\right)\,\left|\,| \hat{{\tilde{u}}}_{k,m}|\le \frac{R}{m^2}\right.\right\},\nonumber
\end{equation}\label{eq:range}
which is compact in $Y_0$. That is, the operator $S$ maps closed convex subsets $Y_0\subset X_0 \subset L_{per}^2((0,T);l^2)$ into compact subsets $Y_{0,c}$ of $Y_0$. The second version of Schauder’s fixed point theorem
implies then that the fixed point equation $u =  S(u)$ has at least one solution. 

We conclude the proof of part (i), with the following detail:
According to the (invertible) transformation (\ref{eq:newx})-(\ref{eq:newp}), breather solutions to system (\ref{eq:transformed}) represent 
breather solutions to the original system (\ref{eq:system}), satisfying the relations (\ref{eq:periodic})-(\ref{eq:limits}). 

Moreover,  by the hypothesis (\ref{eq:nonresonance}) the values of the 
frequency of oscillations $\Omega$ lie above the upper edge  
of the continuous (phonon) spectrum  determined by $2$. Thus the corresponding solutions 
have to be anharmonic which necessitates nonzero nonlinearity. 
Only the latter facilitates amplitude-depending  tuning of the frequency of oscillations so that the breather frequency $\Omega$
lies outside of the phonon spectrum.
Thus, it must  hold that  
$|| N(u)||_{X_0}
\not \equiv 0$, which, by condition {\bf{A}},  is satisfied if and only if  
$u \not \equiv 0$.
That is, the fixed point equation (\ref{eq:compose}) possesses a  nonzero 
solution  and the proof of part (i), is finished.

\textit{(ii): Identification of a ring in $X_0$ for non trivial solutions when \eqref{eq:nonresonance} holds}. Strengthening the condition $\Omega^2>4$ to \eqref{eq:nonresonance}, we are able to identify a ring in $X_0$ for nontrivial solutions. In particular, returning to the estimates \eqref{eq:boundL0},  \eqref{eq:Frechet} and the definition of the map $S$ in \eqref{eq:compose}, we see that from the upper bound \eqref{Ring1} we may derive a contraction regime for the map $S$.  Actually, it follows that if we  require 
$$\frac{ \sqrt{2(1+\cosh(\lambda)}\,\, \overline{K}  
}{\Omega^2-4} R^{\alpha+1}<1,$$
we get that when
\begin{eqnarray}
\label{Ring2}
R\le \left(\frac{\Omega^2-4}{\overline{K}\,\sqrt{2(1+\cosh(\lambda))}}\right)^{\frac{1}{1+\alpha}}:=R_{\mathrm{crit}},
\end{eqnarray} 
there exists a unique solution, namely the trivial one $u=0$. To be consistent with the bound \eqref{eq:As}, and in order to avoid the trivial solution, we need to assume
\begin{eqnarray*}
\label{Ring3}
\Omega^2>4+\overline{K}\,\sqrt{2(1+\cosh(\lambda))},
\end{eqnarray*}
that is, the condition \eqref{eq:nonresonance} on the frequencies holds instead of $\Omega^2>4$. Thus, 
under condition \eqref{eq:nonresonance}, it holds that $R_{\mathrm{crit}}<R_{\mathrm{max}}$ and  non-trivial solutions satisfying \eqref{eq:nonresonance} are located in the ring $\mathbf{R}_E$ of $X_0$ determined by 
\begin{eqnarray}
\label{Ring4}
\mathbf{R}_E=\left\{u\in X_0\;\;:\;\;R_{\mathrm{crit}}\leq ||u||_{X_0}\leq R_{\mathrm{max}}\right\},
\end{eqnarray}
that is the estimates \eqref{eq:As} are satisfied. This concludes the proof of part (ii), and of the whole theorem. \  \ $\Box$
\subsection{Comments on the estimates in $X_0$ and further properties of the solutions}
\label{RemF}
The lower and upper bounds \eqref{eq:As} derived via the fixed-point approach might not be optimal for the estimation of the actual norm of the solutions, as relevant studies concerning DNLS systems (which in some cases may approximate lattice systems described by equation \eqref{eq:system}),  showed \cite{b3,b4,b5}. For example, in the case of the DNLS-type systems, the studies in \cite{b3,b4,b5}  demonstrated that the actual minimum of the power of the solution (excitation threshold), when it exists \cite{Wein99}, is above the relevant non-existence lower bounds similar to $R_{\mathrm{crit}}$ proved herein by contraction mapping-arguments. However, for certain parametric and frequency regimes, the aforementioned studies verified that they can be even sharp estimates of the actual power of the DNLS breathers, particularly in the cases where the excitation threshold (whose emergence depends on the dimension of the lattice \cite{Wein99}) is not existing. In what follows, we discuss features of the estimates, particularly in some limiting cases for their parameters and identify those regimes requiring a coherent change of the parameters.  We also point to  potential ``sharpness'' properties of the estimates in specific parametric regimes.  We conclude by discussing other properties of solutions which follow from the existence Theorem \ref{mainT}. 

\begin{enumerate}
	\item Let us fix $\Omega^2>4$. Considering the limits as $\lambda\rightarrow\infty$ for fixed $\overline{K}$, or as $\overline{K}\rightarrow 0$, for fixed $\lambda$, solely for the upper bound $R_{\mathrm{max}}$,  we observe that 
\begin{eqnarray}
\label{lim1}
\lim_{\lambda\rightarrow\infty} R_{\mathrm{max}}&=&0,\\
\label{lim2}
\lim_{\overline{K}\rightarrow 0} R_{\mathrm{max}}&=&\infty.
\end{eqnarray}
The limit \eqref{lim1} seems to be physically rather relevant in a ``continuous'' regime for the system, than for a discrete regime: In a continuous regime the limit \eqref{lim1} suggests that the more localised the solution becomes, the less is its contribution to its continuous $L^2$-norms. We recall that 
\begin{eqnarray}
\label{norms}
||u||_{l^2}\leq ||u||_{l^2_w},
\end{eqnarray} 
for the sequence spaces, which remains valid as we approximate the continuous case, and is also valid in the continuous limit due to the  weight $w_x=\exp[\lambda |x|]\geq 1$. However, in the discrete regime, the limit \eqref{lim1} is counter-intuitive, since the embedding relation 
\begin{eqnarray}
\label{norms2}
||u||_{l^{\infty}}\leq ||u||_{l^2}\leq ||u||_{l^2_w},
\end{eqnarray}
would imply that the $l^2$-norm of the solution, and consequently, its amplitude,  vanishes.   

On the other hand, the limit \eqref{lim2} seems to be relevant in both the discrete and the continuous regime. Since for $\overline{K}\rightarrow 0$, the system approximates its linear limit, spatially extended  "almost harmonic" modes result instead of localised ones, implying the "unboundness" of the weighted norms.

Both of the above examples implicate a coherent dependence of $\Omega$ on $R$ and  the other parameters. Evidentially 
 the functional dependence of $R$ on all the parameters in the inequality \eqref{Ring1} leads to  the derivation of the estimates.
\item When $\Omega$ satisfies \eqref{eq:nonresonance} instead of only assuming \eqref{eq:nonresonance0}, we may deduce some information concerning the  coherent dependence of the parameters discussed above being physically relevant. 

First, in the limit of large $\lambda$,  condition  \eqref{eq:nonresonance} shows that $\Omega$ becomes arbitrarily large. From the functional expressions of $R_{\mathrm{max}}$ and $R_{\mathrm{crit}}$ we have that 
	\begin{eqnarray}
	\label{Ring5}
	\lim_{(\Omega,\lambda)\rightarrow (\infty,\infty)}R_{\mathrm{crit}}=\lim_{(\Omega,\lambda)\rightarrow (\infty,\infty)}R_{\mathrm{max}}=\infty.	
	\end{eqnarray}
The limit \eqref{Ring5} seems to be in accordance with the physical intuition that, at least for hard interaction potentials, in the limit of arbitrary large frequency,  a type of  ``energy'' of the solution, measured herein in the norm of $X_0$, should become also arbitrarily large. That the weighted norm  grows upon increasing the localisation rate  seems to be consistent with energy concentration phenomena due to enhanced localisation, which may lead accordingly to phenomena of quasi-collapse. Such a phenomenology for the increase of energy was already observed in the formation of a localised state in the $\beta$-FPU model \cite{His}. Note also that in \cite{His} an energy-threshold for the formation of breathers is defined, induced by modulational instability.
\item  The functional forms of the upper and lower bounds in \eqref{Ring5}, which differ only on the $\alpha$-dependent exponents, suggest that they may provide sharp estimates for the actual $X_0$-norm of the localised modes in certain parametric regimes of finite values for the involved parameters. We elaborate further on this, by discussing the features of the assumption \eqref{eq:nonresonance} for small values of $\lambda$ as a potential effective condition for existence of anharmonic oscillations.  Assuming \eqref{eq:nonresonance0}, together with the violation of the condition  \eqref{eq:nonresonance}, implies as $\lambda\rightarrow 0$, that
$$4<\Omega^2 \lesssim 4+2\overline{K}$$
When $\overline{K}\rightarrow 0$, then $\Omega \rightarrow 2^+$. Hence,  as the frequency approaches the upper edge of the phonon band, we do not expect existence of localised modes. Regarding the limit of small $\lambda$, we remark the following: localised solutions on the infinite lattice ${\mathbb{Z}}$ are represented 
by (infinite) square-summable sequences, while exponential decay of the solutions for $|n|\rightarrow \infty$ takes place in the sense of the exponentially weighted 
$l^2$ norm. Notably, for  weight function $w_n\sim 1$, i.e. $\lambda\rightarrow 0$,  our proof establishes the existence  of general localised patterns (e.g. multi-site breathers). 
\item
Since the obtained time-periodic $H^2$ fixed-point-solutions are by Sobolev embeddings $C^1$ in time and since the operator $x \mapsto  W^{\prime}(x)$ maps $C^1$ into itself, one concludes from Eq.\,(\ref{eq:ref1}) that $\ddot{x} \in C^1$   are classical solutions. 
%
\end{enumerate}

\section{Conclusions}
\label{sec4}
We have proved that the Schauder's Fixed Point Theorem can be used to establish the existence  of non-trivial, exponentially localised breather solutions of both even-parity and odd-parity symmetry for  general systems of nonlinearly coupled particles with general interaction potentials, on the infinite lattice ${\mathbb{Z}}$.
The existence problem has been 
reformulated  as a fixed point equation (involving the relevant linear and nonlinear operators associated to the system), on weighted sequence spaces, to directly associate the existence problem with a prescribed rate of spatial localisation, under physically relevant non-resonant conditions, providing in addition, energy bounds for the solutions.
Future extensions may consider systems with other physically important classes of interaction potentials, and importantly (as mentioned in the introductory section), multi-dimensional systems. Another direction may investigate Discrete Nonlinear Schr\"odinger models and different degrees of localisation such as algebraic, which can be relevant for the existence of discrete rational solutions \cite{akhm_AL},\cite{akhm_AL2}. An important task may concern the extension of the method as suitable for the proof of traveling wave solutions \cite{FWattis,Smets,APB}.  Such investigations are in progress and will be documented elsewhere.
\section*{Acknowledgment}
We would like to thank the referee for his/her constructive comments and suggestions. 
\\
\\
\textbf{Authors Declarations}\\
The authors have no conflicts to disclose.\\
\\
\textbf{Data availability}\\
The data that supports the findings of this study are available within the article.

\end{document}